\documentclass{article}
\usepackage{amssymb}
\usepackage{amsfonts}
\usepackage{amsmath}
\usepackage{harvard}

\newtheorem{theorem}{Theorem}

\newtheorem{corollary}[theorem]{Corollary}

\newtheorem{definition}[theorem]{Definition}
\newtheorem{example}[theorem]{Example}

\newtheorem{proposition}[theorem]{Proposition}
\newtheorem{remark}[theorem]{Remark}

\newenvironment{proof}[1][Proof]{\noindent\textbf{#1.} }{\ \rule{0.5em}{0.5em}}
\numberwithin{theorem}{section}
\numberwithin{equation}{section}
\input{tcilatex}

\begin{document}

\title{Canonical Nonlinear Connections in the Multi-Time Hamilton Geometry}
\author{Gheorghe Atanasiu and Mircea Neagu}
\date{}
\maketitle

\begin{abstract}
In this paper we study some geometrical objects (d-tensors, multi-time
semisprays of polymomenta and nonlinear connections) on the dual 1-jet
vector bundle $J^{1\ast }(\mathcal{T},M)\rightarrow \mathcal{T}\times M$.
Some geometrical formulas, which connect the last two geometrical objects,
are also derived. Finally, a canonical nonlinear connection produced by a
Kronecker $h$-regular multi-time Hamiltonian is given.
\end{abstract}



\section{Introduction}

\hspace{4mm} From a geometrical point of view, we point out that the 1-jet
spaces are fundamental ambient mathematical spaces used in the study of
classical and quantum field theories (in their contravariant Lagrangian
approach). For this reason, the differential geometry of these spaces was
intensively studied by many authors (please see, for example, Saunders [15]
or Asanov [1] and references therein). In this direction, it is important to
note that, following the geometrical ideas initially stated by Asanov in
[1], a \textit{multi-time Lagrange contravariant geometry on 1-jet spaces}
(in the sense of distinguished connection, torsions and curvatures) was
recently constructed by Neagu and Udri\c{s}te [12], [14] and published by
Neagu in the book [13]. This geometrical theory is a natural multi-parameter
extension on 1-jet spaces of the already classical \textit{Lagrange
geometrical theory on the tangent bundle} elaborated by Miron and Anastasiei
[10].

From the point of view of physicists, the differential geometry of the dual
1-jet spaces was also studied because the dual 1-jet spaces represent the 
\textit{polymomentum phase spaces} for the \textit{covariant Hamiltonian
formulation of the field theory} (this is a natural \textit{multi-parameter}%
, or \textit{multi-time}, extension of the classical Hamiltonian formalism
from Mechanics). Thus, in order to quantize the covariant Hamiltonian field
theory (this is the final purpose in the framework of quantum field theory),
the \textit{covariant Hamiltonian differential geometry} was developed in
three distinct ways:

\begin{itemize}
\item the \textit{multisymplectic covariant geometry} elaborated by Gotay,
Isenberg, Marsden, Montgomery and their co-workers [6], [7];

\item the \textit{polysymplectic covariant geometry} investigated by
Giachetta, Mangiarotti and Sardanashvily [5];

\item the \textit{De Donder-Weyl covariant Hamiltonian geometry} intensively
studied by Kanatchikov (please see [8] and references therein).
\end{itemize}

These three distinct geometrical-physics variants differ by the multi-time
phase space and the geometrical technics used in study.

Inspired by the Cartan covariant Hamiltonian approach of classical
Mechanics, the studies of Miron [9], Atanasiu [2], [3] and their co-workers
led to the development of the \textit{Hamilton geometry on the cotangent
bundle} exposed in the book [11]. We underline that, via the Legendre
duality of the Hamilton spaces with the Lagrange spaces, it was shown in
[11] that the theory of Hamilton spaces has the same symmetry like the
Lagrange geometry, giving in this way a geometrical framework for the
Hamiltonian theory of Analytical Mechanics.

In such a physical and geometrical context, suggested by the multi-time
framework of the De Donder-Weyl covariant Hamiltonian formulation of
Physical Fields, the aim of this paper is to present some basic geometrical
concepts on dual 1-jet spaces (we refer to distinguished tensors, multi-time
semisprays of polymomenta and nonlinear connections), necessary to the
development of a subsequent \textit{multi-time covariant Hamilton geometry}
(in the sense of \textit{d-connections}, \textit{d-torsions} and \textit{%
d-curvatures} [4]), which to be a natural \textit{multi-parameter}, or 
\textit{poly-momentum}, generalization of the \textit{Hamilton geometry on
the cotangent bundle} [11].

Finally, we would like to point out that the \textit{multi-time Legendre jet
duality }between our subsequent \textit{multi-time covariant Hamilton
geometry} and the already constructed \textit{multi-time contravariant
Lagrange geometry} [13] is a part of our work in progress and represents a
general direction of our future studies.

\section{The dual $1$-jet vector bundle $J^{1\ast }(\mathcal{T},M)$}

\hspace{4mm} We start our geometrical study considering the smooth real
manifolds $\mathcal{T}^{m}$ and $M^{n}$ of dimension $m$, respectively $n$,
whose local coordinates are $(t^{a})_{a=\overline{1,m}}$, respectively $%
(x^{i})_{i=\overline{1,n}}$.

\begin{remark}
i) In this work all geometrical objects and all mappings are considered of
class $C^{\infty }$. This thing is expressed by the words \textbf{%
differentiable} or \textbf{smooth}.

ii) We point out that, throughout this paper, the indices $a,b,c,d,f,g$ run
over the set $\{1,2,\ldots ,m\}$ and the indices $i,j,k,l,r,s$ run over the
set $\{1,2,\ldots ,n\}$.
\end{remark}

Let $(t_{0},x_{0})$ be an arbitrary point of the product manifold $\mathcal{T%
}\times M$ and let $C^{\infty }(\mathcal{T},M)$ be the set all smooth maps
between the manifolds $\mathcal{T}$ and $M$. We define on the space $%
C^{\infty }(\mathcal{T},M)$ the \textit{relation of equivalence} 
\begin{equation*}
\rho \sim _{(t_{0},x_{0})}\sigma \Leftrightarrow \left\{ 
\begin{array}{l}
\rho (t_{0})=\sigma (t_{0})=x_{0}\medskip \\ 
d\rho _{t_{{0}}}=d\sigma _{t_{{0}}}.%
\end{array}%
\right.
\end{equation*}%
Obviously, having two arbitrary smooth maps $\rho ,\sigma \in C^{\infty }(%
\mathcal{T},M)$, the relation of equivalence $\rho \sim _{(t_{{0}},x_{{0}%
})}\sigma $ takes the local form%
\begin{equation*}
\left\{ 
\begin{array}{l}
x^{i}(t_{0}^{b})=y^{i}(t_{0}^{b})=x_{0}^{i}\medskip \\ 
{{\dfrac{\partial x^{i}}{\partial t^{a}}}(t_{0}^{b})={\dfrac{\partial y^{i}}{%
\partial t^{a}}}(t_{0}^{b})},%
\end{array}%
\right.
\end{equation*}%
where $t^{b}(t_{0})=t_{0}^{b},\;x^{i}(x_{0})=x_{0}^{i},\;x^{i}=x^{i}\circ
\rho $ and $y^{i}=x^{i}\circ \sigma $.

The class of equivalence of an element $\rho \in C^{\infty }(\mathcal{T},M)$
is denoted by 
\begin{equation*}
\lbrack \rho ]_{(t_{{0}},x_{{0}})}=\{\sigma \in C^{\infty }(\mathcal{T}%
,M)\;|\;\sigma \sim _{(t_{{0}},x_{{0}})}\rho \}.
\end{equation*}

If we denote by%
\begin{equation*}
J_{t_{{0}},x_{{0}}}^{1}(\mathcal{T},M)=C^{\infty }(\mathcal{T},M)/\sim _{(t_{%
{0}},x_{{0}})},
\end{equation*}%
the quotient space obtained by the factorization of the space $C^{\infty }(%
\mathcal{T},M)$ with respect to the relation of equivalence "$\sim _{(t_{{0}%
},x_{{0}})}$", we can construct the \textit{total space of the jets of order
one}, putting%
\begin{equation*}
J^{1}(\mathcal{T},M)=\bigcup_{(t_{{0}},x_{{0}})\in \mathcal{T}\times M}J_{t_{%
{0}},x_{{0}}}^{1}(\mathcal{T},M).
\end{equation*}

Now, let us organize $J^{1}(\mathcal{T},M)$ like a vector bundle over the
base space $\mathcal{T}\times M$ endowed with the differentiable structure
of a product manifold. For this, we start with an arbitrary smooth map $\rho
\in C^{\infty }(\mathcal{T},M)$, locally given by 
\begin{equation*}
(t^{1},\ldots ,t^{m})\rightarrow (x^{1}(t^{1},\ldots ,t^{m}),\ldots
,x^{n}(t^{1},\ldots ,t^{m})).
\end{equation*}%
Developing the maps $x^{i}$ in Taylor series around the point $%
(t_{0}^{1},\ldots ,t_{0}^{m})\in \mathbb{R}^{m},$ we obtain the local
expressions%
\begin{equation*}
x^{i}(t^{1},\ldots ,t^{m})=x_{0}^{i}+(t^{a}-t_{0}^{a}){\frac{\partial x^{i}}{%
\partial t^{a}}}(t_{0}^{1},\ldots ,t_{0}^{m})+\mathcal{O}(2),
\end{equation*}%
where $(t^{1},\ldots ,t^{m})\in \mathbb{R}^{m}$ is an arbitrary point from a
convenient neighbourhood of the point $(t_{0}^{1},\ldots ,t_{0}^{m})\in 
\mathbb{R}^{m}$, that is $\Vert (t^{1}-t_{0}^{1},\ldots
,t^{m}-t_{0}^{m})\Vert <\varepsilon .$ Considering now the smooth map $%
\overline{\rho }\in C^{\infty }(\mathcal{T},M)$, defined by the set of local
functions%
\begin{equation*}
\overline{x}^{i}(t^{1},\ldots ,t^{m})=x_{0}^{i}+(t^{a}-t_{0}^{a}){\frac{%
\partial x^{i}}{\partial t^{a}}}(t_{0}^{1},\ldots ,t_{0}^{m}),\quad \Vert
(t^{1}-t_{0}^{1},\ldots ,t^{m}-t_{0}^{m})\Vert <\varepsilon ,
\end{equation*}%
we deduce that $\overline{\rho }\sim _{(t_{{0}},x_{{0}})}\rho $. In other
words, the \textit{affine linear approximation} $\overline{\rho }$ of the
map $\rho $ is a very good representative of the class of equivalence $[\rho
]_{(t_{{0}},x_{{0}})}$.

Let $\pi ^{1}:J^{1}(\mathcal{T},M)\rightarrow \mathcal{T}\times M$ be the
canonical projection defined by 
\begin{equation*}
\pi ^{1}([\rho ]_{(t_{{0}},x_{{0}})})=(t_{0},\rho (t_{0})=x_{0}).
\end{equation*}%
It is obvious that the map $\pi ^{1}$ is well defined and surjective. Using
this projection, for each local chart $\mathcal{U}\times V\subset \mathcal{T}%
\times M$, we can define the bijection%
\begin{equation*}
\phi _{\mathcal{U}\times V}:\left( \pi ^{1}\right) ^{-1}(\mathcal{U}\times
V)\rightarrow \mathcal{U}\times V\times \mathbb{R}^{mn},
\end{equation*}%
setting 
\begin{equation*}
\phi _{\mathcal{U}\times V}([\rho ]_{(t_{{0}},x_{{0}})})=\left( t_{0},x_{0},{%
{\frac{\partial x^{i}}{\partial t^{a}}}(t_{0}^{b})}\right) ,\;x_{0}=\rho
(t_{0}).
\end{equation*}

In conclusion, the 1-jet space $J^{1}(\mathcal{T},M)$ can be endowed with a
differentiable structure of dimension $m+n+mn$, such that the maps $\phi _{%
\mathcal{U}\times V}$ to be diffeomorphisms. In this context, the local
coordinates on $J^{1}(\mathcal{T},M)$ are $(t^{a},x^{i},x_{a}^{i})$, where%
\begin{equation*}
\begin{array}{l}
t^{a}([\rho ]_{(t_{{0}},x_{{0}})})=t^{a}(t_{0}),\medskip \\ 
x^{i}([\rho ]_{(t_{{0}},x_{{0}})})=x^{i}(x_{0}),\medskip \\ 
x_{a}^{i}([\rho ]_{(t_{{0}},x_{{0}})})={{\dfrac{\partial x^{i}}{\partial
t^{a}}}(t_{0}^{b}).}%
\end{array}%
\end{equation*}

Using the above coordinates on the 1-jet space $J^{1}(\mathcal{T},M)$, we
get that the projection $\pi ^{1}:J^{1}(\mathcal{T}\times M)\rightarrow 
\mathcal{T}\times M$ has the local expression 
\begin{equation*}
\pi ^{1}(t^{a},x^{i},x_{a}^{i})=(t^{a},x^{i}).
\end{equation*}%
Moreover, the differential map $\pi _{\ast }^{1}$ of the map $\pi ^{1}$ is
locally determined by the Jacobi matrix%
\begin{equation*}
\left( 
\begin{array}{ccc}
\delta _{ab} & 0 & 0 \\ 
0 & \delta _{ij} & 0%
\end{array}%
\right) \in \mathcal{M}_{m+n,m+n+mn}(\mathbb{R}).
\end{equation*}%
Of course, the map $\pi _{\ast }^{1}$ is a surjection (i. e., rank $\pi
_{\ast }^{1}=m+n$) and therefore the map $\pi ^{1}$ is a submersion.
Consequently, we have

\begin{proposition}
The 1-jet space $J^{1}(\mathcal{T},M)$ is a vector bundle over the base
space $\mathcal{T}\times M$, having the fibre type $\mathbb{R}^{mn}$.
\end{proposition}

\begin{remark}
From a physical point of view, the manifold $\mathcal{T}$ can be regarded as
a \textbf{temporal} manifold or, better, a \textbf{multi-time} manifold,
while the manifold $M$ can be regarded as a \textbf{spatial} one. Moreover,
the 1-jet vector bundle%
\begin{equation*}
J^{1}(\mathcal{T},M)\rightarrow \mathcal{T}\times M
\end{equation*}%
can be regarded as a \textbf{bundle of configurations}. This terminology is
justified by the fact that in the particular case $\mathcal{T}=\mathbb{R}$
(i. e., the temporal manifold $\mathcal{T}$ coincides with the usual time
axis represented by the set of real numbers $\mathbb{R}$), we recover the
bundle of configurations which characterizes the classical \textbf{%
non-autonomous}, or \textbf{rheonomic}, \textbf{Mechanics}.
\end{remark}

Taking into account the form of the changes of coordinates on the product
manifold $\mathcal{T}\times M$, we easily deduce

\begin{proposition}
The transformations of coordinates $(t^{a},x^{i},x_{a}^{i})%
\longleftrightarrow (\tilde{t}^{a},\tilde{x}^{i},\tilde{x}_{a}^{i})$ induced
from $\mathcal{T}\times M$ on the 1-jet space $J^{1}(\mathcal{T},M)$ are
given by%
\begin{equation}
\left\{ 
\begin{array}{l}
\tilde{t}^{a}=\tilde{t}^{a}(t^{b})\medskip \\ 
\tilde{x}^{i}=\tilde{x}^{i}(x^{j})\medskip \\ 
\tilde{x}_{a}^{i}={{\dfrac{\partial \tilde{x}^{i}}{\partial x^{j}}}{\dfrac{%
\partial t^{b}}{\partial \tilde{t}^{a}}}x_{b}^{j},}%
\end{array}%
\right.  \label{sch}
\end{equation}%
where $\det (\partial \tilde{t}^{a}/\partial t^{b})\neq 0$ and $\det
(\partial \tilde{x}^{i}/\partial x^{j})\neq 0$.
\end{proposition}

Now, using the general theory of vector bundles (please see [10], for
example), let us consider the \textit{dual 1-jet vector bundle}%
\begin{equation*}
J^{1\ast }(\mathcal{T},M)\rightarrow \mathcal{T}\times M,
\end{equation*}%
whose total space is%
\begin{equation*}
J^{1\ast }(\mathcal{T},M)=\bigcup_{(t_{{0}},x_{{0}})\in \mathcal{T}\times
M}J_{t_{{0}},x_{{0}}}^{1\ast }(\mathcal{T},M),
\end{equation*}%
where 
\begin{equation*}
J_{t_{{0}},x_{{0}}}^{1\ast }(\mathcal{T},M)=\{\omega _{(t_{{0}},x_{{0}%
})}:J_{t_{{0}},x_{{0}}}^{1}(\mathcal{T},M)\rightarrow \mathbb{R}\text{ }|%
\text{ }\omega _{(t_{{0}},x_{{0}})}\text{ is }\mathbb{R}\text{-linear}\},
\end{equation*}%
and which has the fibre type $(\mathbb{R}^{mn})^{\ast }\equiv \mathbb{R}%
^{mn} $. The local coordinates on the dual 1-jet vector bundle $J^{1\ast }(%
\mathcal{T},M)$ are denoted by $(t^{a},x^{i},p_{i}^{a}).$

\begin{remark}
i) In order to simplify the notations, we will use the notations $E=J^{1}(%
\mathcal{T},M)$ and $E^{\ast }=J^{1\ast }(\mathcal{T},M)$.

ii) According to the Kanatchikov's physical terminology [8], which
generalizes the Hamiltonian terminology from Analytical Mechanics, the
coordinates $p_{i}^{a}$ are called \textbf{polymomenta} and the dual 1-jet
space $E^{\ast }$ is called the \textbf{polymomentum phase space}.
\end{remark}

It is easy to see that a transformation of coordinates on the product
manifold $\mathcal{T}\times M$ produces the following results:

\begin{proposition}
The transformations of coordinates $(t^{a},x^{i},p_{i}^{a})%
\longleftrightarrow (\tilde{t}^{a},\tilde{x}^{i},\tilde{p}_{i}^{a})$ induced
from $\mathcal{T}\times M$ on the dual 1-jet space $E^{\ast }$ have the
expressions%
\begin{equation}
\left\{ 
\begin{array}{l}
\tilde{t}^{a}=\tilde{t}^{a}(t^{b})\medskip \\ 
\tilde{x}^{i}=\tilde{x}^{i}(x^{j})\medskip \\ 
\tilde{p}_{i}^{a}={{\dfrac{\partial x^{j}}{\partial \tilde{x}^{i}}}{\dfrac{%
\partial \tilde{t}^{a}}{\partial t^{b}}}}p_{j}^{b}{,}%
\end{array}%
\right.  \label{schp}
\end{equation}%
where $\det (\partial \tilde{t}^{a}/\partial t^{b})\neq 0$ and $\det
(\partial \tilde{x}^{i}/\partial x^{j})\neq 0$.
\end{proposition}

\begin{corollary}
The dual 1-jet space $E^{\ast }$ is an orientable manifold having the
dimension $m+n+mn.$
\end{corollary}

Now, doing a transformation of coordinates (\ref{schp}) on $E^{\ast }$, we
obtain

\begin{proposition}
The elements of the local natural basis 
\begin{equation*}
{\left\{ {\frac{\partial }{\partial t^{a}}},{\frac{\partial }{\partial x^{i}}%
},{\frac{\partial }{\partial p_{i}^{a}}}\right\} }
\end{equation*}%
of the Lie algebra of vector fields $\mathcal{X}(E^{\ast })$ transform by
the rules%
\begin{equation}
\begin{array}{l}
{{\dfrac{\partial }{\partial t^{a}}}={\dfrac{\partial \tilde{t}^{b}}{%
\partial t^{a}}}{\dfrac{\partial }{\partial \tilde{t}^{b}}}+{\dfrac{\partial 
\tilde{p}_{j}^{b}}{\partial t^{a}}}{\dfrac{\partial }{\partial \tilde{p}%
_{j}^{b}}}},\medskip \\ 
{{\dfrac{\partial }{\partial x^{i}}}={\dfrac{\partial \tilde{x}^{j}}{%
\partial x^{i}}}{\dfrac{\partial }{\partial \tilde{x}^{j}}}+{\dfrac{\partial 
\tilde{p}_{j}^{b}}{\partial x^{i}}}{\dfrac{\partial }{\partial \tilde{p}%
_{j}^{b}}}},\medskip \\ 
{{\dfrac{\partial }{\partial p_{i}^{a}}}={\dfrac{\partial x^{i}}{\partial 
\tilde{x}^{j}}}{\dfrac{\partial \tilde{t}^{b}}{\partial t^{a}}}{\dfrac{%
\partial }{\partial \tilde{p}_{j}^{b}}}.}%
\end{array}
\label{schc}
\end{equation}
\end{proposition}

\begin{proposition}
The elements of the local natural cobasis $\{dt^{a},dx^{i},dp_{i}^{a}\}$ of
the Lie algebra of covector fields $\mathcal{X}^{\ast }(E^{\ast })$
transform by the rules%
\begin{equation}
\begin{array}{l}
{dt^{a}={\dfrac{\partial t^{a}}{\partial \tilde{t}^{b}}}d\tilde{t}^{b}}%
,\medskip \\ 
{dx^{i}={\dfrac{\partial x^{i}}{\partial \tilde{x}^{j}}}d\tilde{x}^{j}}%
,\medskip \\ 
{dp_{i}^{a}={\dfrac{\partial p_{i}^{a}}{\partial \tilde{t}^{b}}}d\tilde{t}%
^{b}+{\dfrac{\partial p_{i}^{a}}{\partial \tilde{x}^{j}}}d\tilde{x}^{j}+{%
\dfrac{\partial \tilde{x}^j}{\partial x^{i}}}{\dfrac{\partial t^{a}}{%
\partial \tilde{t}^b}}d\tilde{p}_j^b.}%
\end{array}
\label{schf}
\end{equation}
\end{proposition}

\begin{remark}
Let us remark that, in the particular case $\mathcal{T}=\mathbb{R},$ we find
the \textbf{momentum phase space} $J^{1\ast }(\mathbb{R},M)\equiv \mathbb{R}%
\times T^{\ast }M$ (we have a punctual identification)$,$ where $T^{\ast }M$
is the cotangent bundle. This particular momentum phase space $\mathbb{R}%
\times T^{\ast }M$ is regarded as a vector bundle over the product of
manifolds $\mathbb{R}\times M.$ Its coordinates are denoted by $%
(t,x^{i},p_{i})$ and the corresponding transformations group (\ref{schp})
becomes%
\begin{equation}
\left\{ 
\begin{array}{l}
\tilde{t}=\tilde{t}(t)\medskip \\ 
\tilde{x}^{i}=\tilde{x}^{i}(x^{j})\medskip \\ 
\tilde{p}_{i}={{\dfrac{\partial x^{j}}{\partial \tilde{x}^{i}}}{\dfrac{d%
\tilde{t}}{dt}}}p_{j}{.}%
\end{array}%
\right.  \label{schp1}
\end{equation}

It is important to note that the group of transformations (\ref{schp1})
emphasizes the \textbf{relativistic} character played by the usual time $t$,
and it is different by the group of transformations 
\begin{equation}
\left\{ 
\begin{array}{l}
\tilde{t}=t\medskip \\ 
\tilde{x}^{i}=\tilde{x}^{i}\left( x^{j}\right) ,\medskip \\ 
\tilde{p}_{i}=\dfrac{\partial x^{j}}{\partial \tilde{x}^{i}}p_{j}%
\end{array}%
\right.  \label{tr-group-absolut-t}
\end{equation}%
of the trivial bundle $\mathbb{R}\times T^{\ast }M\rightarrow T^{\ast }M$
used in the \textbf{non-autonomous Hamilton geometry }for the study of the
metrical structure%
\begin{equation*}
g^{ij}\left( t,x,p\right) =\frac{1}{2}\frac{\partial ^{2}H}{\partial
p_{i}\partial p_{j}},
\end{equation*}%
where $H:\mathbb{R}\times T^{\ast }M\rightarrow \mathbb{R}$ is a Hamiltonian
function.
\end{remark}

\section{d-Tensors, multi-time semisprays of polymomenta and nonlinear
connections}

\hspace{4mm} It is well known the importance of tensors in the development
of a fertile geometry on a vector bundle. Following the geometrical ideas
developed in the books [10] and [11], in our study upon the geometry of the
dual 1-jet bundle $E^{\ast }$ a central role is played by the \textit{%
distinguished tensors} or, briefly, \textit{d-tensors}.

\begin{definition}
A geometrical object $T=\left( T_{bj(c)(l)\ldots }^{ai(k)(d)\ldots }\right) $
on the dual 1-jet vector bundle $E^{\ast }$, whose local components, with
respect to a transformation of coordinates (\ref{schp}) on $E^{\ast }$,
transform by the rules%
\begin{equation*}
T_{bj(c)(l)\ldots }^{ai(k)(d)\ldots }=\tilde{T}_{fq(g)(s)\ldots
}^{ep(r)(h)\ldots }{\frac{\partial t^{a}}{\partial \tilde{t}^{e}}}{\frac{%
\partial x^{i}}{\partial \tilde{x}^{p}}}\left( {\frac{\partial x^{k}}{%
\partial \tilde{x}^{r}}}{\frac{\partial \tilde{t}^{g}}{\partial t^{c}}}%
\right) {\frac{\partial \tilde{t}^{f}}{\partial t^{b}}}{\frac{\partial 
\tilde{x}^{q}}{\partial x^{j}}}\left( {\frac{\partial \tilde{x}^{s}}{%
\partial x^{l}}}{\frac{\partial t^{d}}{\partial \tilde{t}^{h}}}\right)
\ldots \;,
\end{equation*}%
is called a \textbf{d-tensor} or \textbf{distinguished tensor field} on the
dual 1-jet space $E^{\ast }$.
\end{definition}

\begin{remark}
The utilization between parentheses of certain indices of the local
components $T_{bj(c)(l)\ldots }^{ai(k)(d)\ldots }$ is necessary for clearer
future contractions\textit{.} For the moment, we point out only that ${%
\QATOP{(k)}{(c)}}$ or ${\QATOP{(d)}{(l)}}$ behaves like a single \textbf{%
double index}.
\end{remark}

\begin{example}
i) If $H:E^{\ast }\rightarrow \mathbb{R}$ is a Hamiltonian function
depending on the polymomenta $p_{i}^{a}$, then the local components%
\begin{equation*}
G_{(a)(b)}^{(i)(j)}={\frac{1}{2}}{\frac{\partial ^{2}H}{\partial
p_{i}^{a}\partial p_{j}^{b}}}
\end{equation*}%
represent a d-tensor field $\mathbb{G}=\left( G_{(a)(b)}^{(i)(j)}\right) $
on the dual 1-jet space $E^{\ast }$, which is called the \textbf{fundamental
vertical metrical d-tensor associated to the Hamiltonian function of
polymomenta} $H.$ This is because if $\mathcal{T}=\mathbb{R}$ and $H$ is a
regular Hamiltonian function, then the d-tensor field $\mathbb{G}$ can be
regarded as the fundamental metrical d-tensor $g^{ij}(t,x,p)$ from the
theory of rheonomic Hamilton spaces.

ii) Let us consider the d-tensor $\mathbb{C}^{\ast }=\left( \mathbb{C}%
_{(i)}^{(a)}\right) $, where $\mathbb{C}_{(i)}^{(a)}=p_{i}^{a}$. The
distinguished tensor $\mathbb{C}^{\ast }$ is called the \textbf{%
Liou\-ville-Ha\-mil\-ton d-tensor field of polymomenta} on the dual 1-jet
space $E^{\ast }$. Remark that, for the particular case $\mathcal{T}=\mathbb{%
R}$, we recover the classical Liouville-Hamilton vector field 
\begin{equation*}
\mathbb{C}^{\ast }\mathbb{=}p_{i}{{\frac{\partial }{\partial p_{i}}}}
\end{equation*}%
on the cotangent bundle $T^{\ast }M$, which is used in the Hamilton geometry
[11].

iii) Let $h_{ab}(t)$ be a semi-Riemannian metric on the temporal manifold $%
\mathcal{T}$. The geometrical object $\mathbb{L}=\left(
L_{(j)ab}^{(c)}\right) $, where 
\begin{equation*}
L_{(j)ab}^{(c)}=h_{ab}p_{j}^{c},
\end{equation*}%
is a d-tensor field on $E^{\ast }$, which is called the \textbf{polymomentum
Liou\-ville-Ha\-mil\-ton d-tensor field associated to the metric }$h_{ab}(t)$%
.

iv) Using the preceding metric $h_{ab}(t)$, we can construct the d-tensor
field $\mathbb{J}=\left( J_{(a)bj}^{(i)}\right) $, where 
\begin{equation*}
J_{(a)bj}^{(i)}=h_{ab}\delta _{j}^{i}.
\end{equation*}%
The distinguished tensor $\mathbb{J}$ is called the \textbf{d-tensor of }$h$%
\textbf{-normalization }\textit{on the dual 1-jet vector bundle }$E^{\ast }$.
\end{example}

It is obvious that any d-tensor field on $E^{\ast }$ is a tensor field on $%
E^{\ast }$. Conversely, this statement is not true. As examples, we
construct two tensors on $E^{\ast }$, which are not distinguished tensors on 
$E^{\ast }$.

\begin{definition}
A global tensor $\underset{1}{G}$ on $E^{\ast },$ locally expressed by%
\begin{equation*}
\underset{1}{G}=p_{i}^{a}dx^{i}\otimes {\frac{\partial }{\partial t^{a}}}-2%
\underset{1}{G}\text{{}}_{(j)i}^{(b)}dx^{i}\otimes {\frac{\partial }{%
\partial p_{j}^{b}}},
\end{equation*}%
is called a \textbf{temporal semispray} on the dual 1-jet vector bundle $%
E^{\ast }$.
\end{definition}

Taking into account that the temporal semispray $\underset{1}{G}$ is a
global tensor on $E^{\ast }$, by a direct calculation, we obtain

\begin{proposition}
i) With respect to a transformation of coordinates (\ref{schp}), the
components $\underset{1}{G}${}$_{(j)i}^{(b)}$ of the global tensor $\underset%
{1}{G}$ transform by the rules%
\begin{equation}
2\underset{1}{\widetilde{G}}\text{{}}_{(k)r}^{(c)}=2\underset{1}{G}\text{{}}%
_{(j)i}^{(b)}{\frac{\partial \tilde{t}^{c}}{\partial t^{b}}}{\frac{\partial
x^{i}}{\partial \tilde{x}^{r}}}{\frac{\partial x^{j}}{\partial \tilde{x}^{k}}%
}-{\frac{\partial x^{i}}{\partial \tilde{x}^{r}}}{\frac{\partial \tilde{p}%
_{k}^{c}}{\partial t^{a}}}p_{i}^{a}.  \label{tspr}
\end{equation}%
In other words, the temporal semispray $\underset{1}{G}$ is not a d-tensor.

ii) Conversely, to give a temporal semispray on $E^{\ast }$ is equivalent to
give a set of local functions $\underset{1}{G}=\left( \underset{1}{G}%
{}_{(j)i}^{(b)}\right) $ which transform by the rules (\ref{tspr}).
\end{proposition}

Moreover, if we start with 
\begin{equation*}
\underset{1}{G}=\underset{1}{G}\text{{}}_{i}^{a}dx^{i}\otimes {\frac{%
\partial }{\partial t^{a}}}-2\underset{1}{G}{}_{(j)i}^{(b)}dx^{i}\otimes {%
\frac{\partial }{\partial p_{j}^{b}}}
\end{equation*}%
a global tensor on $E^{\ast }$ and $h_{ab}(t)$ an arbitrary semi-Riemannian
metric on the temporal manifold $\mathcal{T}$, then we easily find the
following result:

\begin{proposition}
The tensor $\underset{1}{G}$ is a temporal semispray on $E^{\ast }$ if and
only if 
\begin{equation*}
J_{(a)bj}^{(i)}\underset{1}{G}\text{{}}_{i}^{c}=L_{(j)ab}^{(c)},
\end{equation*}%
where $\mathbb{J}$ is the d-tensor of $h$-normalization and $\mathbb{L}$ is
the multi-time Liou\-ville-Ha\-mil\-ton d-tensor field associated to the
metric $h_{ab}(t)$.
\end{proposition}

\begin{example}
If $\varkappa _{bc}^{a}(t)$ are the Christoffel symbols of a semi-Riemannian
me\-tric $h_{ab}(t)$ of the temporal manifold $\mathcal{T}$, then the local
components 
\begin{equation}
\underset{1}{\overset{0}{G}}\text{{}}_{(j)k}^{(a)}=\frac{1}{2}\varkappa
_{bc}^{a}p_{j}^{b}p_{k}^{c}  \label{temp-semispray-assoc-metric}
\end{equation}%
represent a temporal semispray $\underset{1}{\overset{0}{G}}$ on the dual
1-jet vector bundle $E^{\ast }$.
\end{example}

\begin{definition}
The temporal semispray $\underset{1}{\overset{0}{G}}$ given by (\ref%
{temp-semispray-assoc-metric}) is called the \textbf{canonical temporal
semispray associated to metric }$h_{ab}(t)$\textit{.}
\end{definition}

A second example of tensor on the dual 1-jet space $E^{\ast },$ which is not
a distinguished tensor, is offered by

\begin{definition}
A global tensor $\underset{2}{G}$ on $E^{\ast },$ locally expressed by%
\begin{equation*}
\underset{2}{G}=\delta _{i}^{j}dx^{i}\otimes {\frac{\partial }{\partial x^{j}%
}}-2\underset{2}{G}\text{{}}_{(j)i}^{(b)}dx^{i}\otimes {\frac{\partial }{%
\partial p_{j}^{b}},}
\end{equation*}%
is called a \textbf{spatial semispray} on the dual 1-jet vector bundle $%
E^{\ast }$.
\end{definition}

As in the case of a temporal semispray, we can prove without difficulties
the following statements:

\begin{proposition}
i) To give a spatial semispray on $E^{\ast }$ is equivalent to give a set of
local functions $\underset{2}{G}=\left( \underset{2}{G}\text{{}}%
_{(j)i}^{(b)}\right) $ which transform by the rules%
\begin{equation}
2\underset{2}{\widetilde{G}}\text{{}}_{(s)k}^{(d)}=2\underset{2}{G}\text{{}}%
_{(j)i}^{(b)}{\frac{\partial \tilde{t}^{d}}{\partial t^{b}}}{\frac{\partial
x^{i}}{\partial \tilde{x}^{k}}}{\frac{\partial x^{j}}{\partial \tilde{x}^{s}}%
}-{\frac{\partial x^{i}}{\partial \tilde{x}^{k}}}{\frac{\partial \tilde{p}%
_{s}^{d}}{\partial x^{i}}}.  \label{sspr}
\end{equation}

ii) A global tensor on $E^{\ast }$, defined by 
\begin{equation*}
\underset{2}{G}=\underset{2}{G}\text{{}}_{i}^{j}dx^{i}\otimes {\frac{%
\partial }{\partial x^{j}}}-2\underset{2}{G}\text{{}}_{(j)i}^{(b)}dx^{i}%
\otimes {\frac{\partial }{\partial p_{j}^{b}},}
\end{equation*}%
is a spatial semispray on $E^{\ast }$ if and only if 
\begin{equation*}
J_{(a)bj}^{(i)}\underset{2}{G}\text{{}}_{i}^{k}=J_{(a)bj}^{(k)},
\end{equation*}%
where $\mathbb{J}$ is the d-tensor of $h$-normalization for an arbitrary
semi-Riemannian temporal metric $h_{ab}(t)$.
\end{proposition}

\begin{example}
If $\gamma _{jk}^{i}(x)$ are the Christoffel symbols of a semi-Riemannian
me\-tric $\varphi _{ij}(x)$ of the spatial manifold $M$, then the local
components 
\begin{equation}
\underset{2}{\overset{0}{G}}\text{{}}_{(j)k}^{(b)}=-\frac{1}{2}\gamma
_{jk}^{i}p_{i}^{b}  \label{spatial-semispr-assoc-metric}
\end{equation}%
define a spatial semispray $\underset{2}{\overset{0}{G}}$ on the dual 1-jet
space $E^{\ast }$.
\end{example}

\begin{definition}
The spatial semispray $\underset{2}{\overset{0}{G}}$ given by (\ref%
{spatial-semispr-assoc-metric}) is called the \textbf{canonical spatial
semispray associated to the metric }$\varphi _{ij}(x)$\textit{.}
\end{definition}

\begin{remark}
\label{Rem} It is obvious that the difference between two temporal (spatial,
respectively) semisprays is a d-tensor.
\end{remark}

Using the Remark \ref{Rem} and the preceding notations, we easily deduce

\begin{theorem}
\label{th2} Let $(\mathcal{T},h)$ and $(M,\varphi )$ be two semi-Riemannian
manifolds and let $\underset{1}{G}=\left( \underset{1}{G}{}_{(i)j}^{(b)}%
\right) $ ($\underset{2}{G}=\left( \underset{2}{G}\text{{}}%
_{(i)j}^{(b)}\right) $, respectively) be an arbitrary temporal (spatial,
respectively) semispray on the dual 1-jet space $E^{\ast }$. In this
context, the following equalities are true:%
\begin{equation*}
\begin{array}{ll}
\underset{1}{G}{}_{(i)j}^{(b)}=\dfrac{1}{2}\varkappa
_{cd}^{b}p_{i}^{c}p_{j}^{d}+\underset{1}{T}\text{{}}_{(i)j}^{(b)}, & 
\underset{2}{G}\text{{}}_{(i)j}^{(b)}=-\dfrac{1}{2}\gamma _{ij}^{k}p_{k}^{b}+%
\underset{2}{T}\text{{}}_{(i)j}^{(b)},%
\end{array}%
\end{equation*}%
where $\underset{1}{T}${}$_{(i)j}^{(b)},\;\underset{2}{T}${}$_{(i)j}^{(b)}$
are unique d-tensors with the preceding properties.
\end{theorem}

\begin{definition}
A pair $G=\left( \underset{1}{G},\underset{2}{G}\right) $ consisting of a
temporal semispray $\underset{1}{G}$ and a spatial one $\underset{2}{G}$ is
called a \textbf{multi-time semispray of polymomenta} on the dual 1-jet
space $E^{\ast }$.
\end{definition}

\begin{remark}
The Theorem \ref{th2} emphasizes the central role played by the canonical
semispray of polymomenta $\overset{0}{G}=\left( \underset{1}{\overset{0}{G}},%
\underset{2}{\overset{0}{G}}\right) $, associated to a pair of
semi-Riemannian metrics $(h_{ab}(t),\varphi _{ij}(x))$, in the description
of an arbitrary multi-time semispray of polymomenta $G=\left( \underset{1}{G}%
,\underset{2}{G}\right) $ on the dual 1-jet space $E^{\ast }$.
\end{remark}

\begin{definition}
A pair of local functions $N=\left( \underset{1}{N}\text{{}}_{(k)a}^{(c)},%
\underset{2}{N}\text{{}}_{(k)i}^{(c)}\right) $ on $E^{\ast }$, which
transform by the rules%
\begin{equation}
\begin{array}{l}
\underset{1}{\widetilde{N}}\text{{}}{_{(j)d}^{(b)}=}\underset{1}{N}\text{{}}%
_{(k)a}^{(c)}\dfrac{\partial \tilde{t}^{b}}{\partial t^{c}}{{\dfrac{\partial
x^{k}}{\partial \tilde{x}^{j}}}{\dfrac{\partial t^{a}}{\partial \tilde{t}^{d}%
}}-\dfrac{\partial t^{a}}{\partial \tilde{t}^{d}}{\dfrac{\partial \tilde{p}%
_{j}^{b}}{\partial t^{a}}}},\medskip \\ 
\underset{2}{\widetilde{N}}\text{{}}{_{(j)r}^{(b)}=\underset{2}{N}%
{}_{(k)i}^{(c)}{\dfrac{\partial \tilde{t}^{b}}{\partial t^{c}}}{\dfrac{%
\partial x^{k}}{\partial \tilde{x}^{j}}}}\dfrac{\partial x^{i}}{\partial 
\tilde{x}^{r}}{-\dfrac{\partial x^{i}}{\partial \tilde{x}^{r}}{\dfrac{%
\partial \tilde{p}_{j}^{b}}{\partial x^{i}}}},%
\end{array}
\label{schco}
\end{equation}%
is called a \textbf{nonlinear connection} on the dual 1-jet bundle $E^{\ast
} $.
\end{definition}

\begin{remark}
The nonlinear connections are very important in the study of the
differential geometry of the dual 1-jet space $E^{\ast }$ because they
produce the \textbf{adapted distinguished 1-forms}%
\begin{equation*}
\delta p_{i}^{a}=dp_{i}^{a}+\underset{1}{N}\text{{}}_{(i)b}^{(a)}dt^{b}+{%
\underset{2}{N}{}_{(i)j}^{(a)}}dx^{j},
\end{equation*}%
which are necessary for the adapted local description of the geometrical
objects involved in study, such as the \textbf{d-connections}, the \textbf{%
d-torsions} or the \textbf{d-curvatures}. For more details, please see the
paper [4].
\end{remark}

Now, let us expose the connection between the notions of multi-time
semispray of polymomenta and nonlinear connection on the dual 1-jet space $%
E^{\ast }$. Thus, in our context, using the transformation rules (\ref{tspr}%
), (\ref{sspr}) and (\ref{schco}) of the geometrical objects taken in study,
we can easily prove the following statements:

\begin{proposition}
\label{scospr} i) If $\underset{1}{G}${}$_{(j)k}^{(a)}$ are the components
of a temporal semispray $\underset{1}{G}$ on $E^{\ast }$ and $\varphi
_{ij}(x)$ is a semi-Riemannian metric on the spatial manifold $M$, then the
local components%
\begin{equation*}
\underset{1}{N}\text{{}}_{(r)b}^{(a)}=\varphi ^{jk}{\frac{\partial \underset{%
1}{G}{}_{(j)k}^{(a)}}{\partial p_{i}^{b}}}\varphi _{ir}
\end{equation*}%
represent the temporal components of a nonlinear connection $N_{G}$ on $%
E^{\ast }$.

ii) Conversely, if $\underset{1}{N}${}$_{(i)b}^{(a)}$ are the temporal
components of a nonlinear connection $N$ on $E^{\ast }$, then the local
components%
\begin{equation*}
\underset{1}{G}\text{{}}_{(i)j}^{(a)}=\frac{1}{2}\underset{1}{N}%
{}_{(i)b}^{(a)}p_{j}^{b}
\end{equation*}%
represent a temporal semispray $\underset{1}{G}${}$_{N}$ on $E^{\ast }$.
\end{proposition}

\begin{proposition}
\label{tcospr} i) If $\underset{2}{G}$ $\!\!_{(j)i}^{(b)}$ are the
components of a spatial semispray $\underset{2}{G}$ on $E^{\ast }$, then the
local components%
\begin{equation*}
\underset{2}{N}\text{{}}_{(j)i}^{(b)}=2\underset{2}{G}\text{{}}_{(j)i}^{(b)}
\end{equation*}%
represent the spatial components of a nonlinear connection $N_{G}$ on $%
E^{\ast }$.

ii) Conversely, if $\underset{2}{N}${}$_{(j)i}^{(b)}$ are the spatial
components of a nonlinear connection $N$ on $E^{\ast }$, then the local
functions 
\begin{equation*}
\underset{2}{G}\text{{}}_{(j)i}^{(b)}{={\frac{1}{2}}}\underset{2}{N}\text{{}}%
_{(j)i}^{(b)}
\end{equation*}%
represent a spatial semispray $\underset{2}{G}${}$_{N}$ on $E^{\ast }$.
\end{proposition}

\begin{remark}
The Propositions \ref{scospr} and \ref{tcospr} emphasize that a multi-time
semispray of polymomenta $G=\left( \underset{1}{G},\underset{2}{G}\right) $
on the dual 1-jet space $E^{\ast }$ naturally induces a nonlinear connection 
$N_{G}$ on $E^{\ast }$ and vice-versa, $N$ induces $G_{N}$.
\end{remark}

\begin{definition}
The nonlinear connection $N_{G}$ on the dual 1-jet space $E^{\ast }$ is
called the \textbf{canonical nonlinear connection associated to the
multi-time semispray of polymomenta} $G=\left( \underset{1}{G},\underset{2}{G%
}\right) $ and vice-versa.
\end{definition}

\begin{corollary}
\label{NLC0} The canonical nonlinear connection $\overset{0}{N}=\left( 
\underset{1}{\overset{0}{N}}\text{{}}_{(i)b}^{(a)},\underset{2}{\overset{0}{N%
}}\text{{}}_{(i)j}^{(a)}\right) $ produced by the canonical multi-time
semispray of polymomenta $\overset{0}{G}=\left( \underset{1}{\overset{0}{G}},%
\underset{2}{\overset{0}{G}}\right) $ associated to the pair of
semi-Riemannian metrics $(h_{ab}(t),\varphi _{ij}(x))$ has the local
components%
\begin{equation*}
\begin{array}{ccc}
\underset{1}{\overset{0}{N}}\text{{}}_{(i)b}^{(a)}=\varkappa
_{cb}^{a}p_{i}^{c} & \text{and} & \underset{2}{\overset{0}{N}}\text{{}}%
_{(i)j}^{(a)}=-\gamma _{ij}^{k}p_{k}^{a}.%
\end{array}%
\end{equation*}
\end{corollary}

\section{Kronecker $h$-regularity. Canonical nonlinear connections}

\hspace{4mm} Let us consider a smooth multi-time Hamiltonian function $%
H:E^{\ast }\rightarrow \mathbb{R}$, locally expressed by%
\begin{equation*}
E^{\ast }\ni (t^{a},x^{i},p_{i}^{a})\rightarrow H(t^{a},x^{i},p_{i}^{a})\in 
\mathbb{R},
\end{equation*}%
whose \textit{fundamental vertical metrical d-tensor} is defined by%
\begin{equation*}
G_{(a)(b)}^{(i)(j)}={\frac{1}{2}}{\frac{\partial ^{2}H}{\partial
p_{i}^{a}\partial p_{j}^{b}}}.
\end{equation*}

In the sequel, let us fix $h=(h_{ab}(t^{c})),$ a semi-Riemannian metric on
the temporal manifold $\mathcal{T}$, together with a d-tensor $%
g^{ij}(t^{c},x^{k},p_{k}^{c})$ on the dual 1-jet space $E^{\ast }$, which is
symmetric, has the rank $n=\dim M$ and a constant signature.

\begin{definition}
A multi-time Hamiltonian function $H:E^{\ast }\rightarrow \mathbb{R},$
having the fundamental vertical metrical d-tensor of the form%
\begin{equation*}
G_{(a)(b)}^{(i)(j)}(t^{c},x^{k},p_{k}^{c})={\frac{1}{2}}{\frac{\partial ^{2}H%
}{\partial p_{i}^{a}\partial p_{j}^{b}}}%
=h_{ab}(t^{c})g^{ij}(t^{c},x^{k},p_{k}^{c}),
\end{equation*}%
is called a \textbf{Kronecker }$h$\textbf{-regular multi-time Hamiltonian
function}\textit{.}
\end{definition}

In this context, we can introduce the following geometrical concept:

\begin{definition}
A pair $MH_{m}^{n}=(E^{\ast }=J^{1\ast }(\mathcal{T},M),H),$ where $m=\dim 
\mathcal{T}$ and $n=\dim M,$ consisting of the dual 1-jet space and a
Kronecker $h$-regular multi-time Hamiltonian function $H:E^{\ast
}\rightarrow \mathbb{R},$ is called a \textbf{multi-time Hamilton space}.
\end{definition}

\begin{remark}
i) In the particular case $(\mathcal{T},h)=(\mathbb{R},\delta ),$ a
multi-time Hamilton space will be called a \textbf{relativistic rheonomic
Hamilton space}\textit{\ and it will be denoted by }$RH^{n}=(J^{1\ast }(%
\mathbb{R},M),H)$.

ii) If the temporal manifold $(\mathcal{T},h)$ is 1-dimensional, then, via a
temporal reparametrization, we have $J^{1\ast }(\mathcal{T},M)\equiv
J^{1\ast }(\mathbb{R},M)$. In other words, a multi-time Hamilton space
having $\dim \mathcal{T}=1$ is a reparametrized relativistic rheonomic
Hamilton space.
\end{remark}

\begin{example}
Let us consider the following Kronecker $h$-regular multi-time Hamiltonian
function $H_{1}:E^{\ast }\rightarrow \mathbb{R},$ defined by%
\begin{equation}
H_{1}=\frac{1}{mc}h_{ab}(t)\varphi ^{ij}(x)p_{i}^{a}p_{j}^{b},  \label{G}
\end{equation}%
where $h_{ab}(t)$ ($\varphi _{ij}(x)$, respectively) is a semi-Riemannian
metric on the temporal (spatial, respectively) manifold $\mathcal{T}$ ($M$,
respectively) having the physical meaning of \textbf{gravitational potentials%
}, and $m$ and $c$ are the known constants from Physics representing the 
\textbf{mass of the test body} and the \textbf{speed of light}. Then, the
multi-time Hamilton space 
\begin{equation*}
\mathcal{G}MH_{m}^{n}=(E^{\ast },H_{1})
\end{equation*}%
defined by the multi-time Hamiltonian function (\ref{G}) is called the 
\textbf{multi-time Hamilton space of gravitational field}. This is because,
for $(\mathcal{T},h)=(\mathbb{R},\delta )$, we recover the classical
Hamilton space of gravitational field from the book [11].
\end{example}

\begin{example}
Using preceding notations, let us consider the Kronecker $h$-regular
multi-time Hamiltonian function $H_{2}:E^{\ast }\rightarrow \mathbb{R},$
defined by%
\begin{equation}
H_{2}=\frac{1}{mc}h_{ab}(t)\varphi ^{ij}(x)p_{i}^{a}p_{j}^{b}-\frac{2e}{%
mc^{2}}A_{(a)}^{(i)}(x)p_{i}^{a}+\frac{e^{2}}{mc^{3}}F(t,x),  \label{ED}
\end{equation}%
where $A_{(a)}^{(i)}(x)$ is a d-tensor on $E^{\ast }$ having the physical
meaning of \textbf{potential d-tensor of an electromagnetic field}, $e$ is
the \textbf{charge of the test body} and the function $F(t,x)$ is given by%
\begin{equation*}
F(t,x)=h^{ab}(t)\varphi _{ij}(x)A_{(a)}^{(i)}(x)A_{(b)}^{(j)}(x).
\end{equation*}%
Then, the multi-time Hamilton space 
\begin{equation*}
\mathcal{ED}MH_{m}^{n}=(E^{\ast },H_{2})
\end{equation*}%
defined by the multi-time Hamiltonian function (\ref{ED}) is called the 
\textbf{autonomous multi-time Hamilton space of electrodynamics}. This is
because, in the particular case $(\mathcal{T},h)=(\mathbb{R},\delta )$, we
recover the classical Hamilton space of electrodynamics studied in the book
[11]. The non-dynamical character (the independence of the temporal
coordinates $t^{c}$) of the spatial gravitational potentials $\varphi
_{ij}(x)$ motivated us to use the term \textbf{"autonomous"}.
\end{example}

\begin{example}
More general, if we take on $E^{\ast }$ a symmetric d-tensor field $%
g_{ij}(t,x)$ having the rank $n$ and a constant signature, we can define the
Kronecker $h$-regular multi-time Hamiltonian function $H_{3}:E^{\ast
}\rightarrow \mathbb{R},$ setting%
\begin{equation}
H_{3}=h_{ab}(t)g^{ij}(t,x)p_{i}^{a}p_{j}^{b}+U_{(a)}^{(i)}(t,x)p_{i}^{a}+%
\mathcal{F}(t,x),  \label{NED}
\end{equation}%
where $U_{(a)}^{(i)}(t,x)$ is a d-tensor field on $E^{\ast }$ and $\mathcal{F%
}(t,x)$ is a function on $E^{\ast }$. Then, the multi-time Hamilton space 
\begin{equation*}
\mathcal{NED}MH_{m}^{n}=(E^{\ast },H_{3})
\end{equation*}%
defined by the multi-time Hamiltonian function (\ref{NED}) is called the 
\textbf{non-a\-u\-to\-no\-mous multi-time Hamilton space of electrodynamics}%
. The dynamical character (the dependence of the temporal coordinates $t^{c}$%
) of the gravitational potentials $g_{ij}(t,x)$ motivated us to use the word 
\textbf{"non-autonomous".}
\end{example}

An important role and, at the same time, an obstruction for the subsequent
development of a geometrical theory of the multi-time Hamilton spaces, is
re\-pre\-sen\-ted by

\begin{theorem}[of characterization of multi-time Hamilton spaces]
\label{thchar} If\linebreak we have $m=\dim \mathcal{T}\geq 2$, then the
following statements are equivalent:

(i) $H$ is a Kronecker $h$-regular multi-time Hamiltonian function on $%
E^{\ast }$.

(ii) The multi-time Hamiltonian function $H$ reduces to a multi-time
Hamiltonian function of non-autonomous electrodynamic kind, that is we have%
\begin{equation}
H=h_{ab}(t)g^{ij}(t,x)p_{i}^{a}p_{j}^{b}+U_{(a)}^{(i)}(t,x)p_{i}^{a}+%
\mathcal{F}(t,x).  \label{NEDTH}
\end{equation}
\end{theorem}

\begin{proof}
(ii) $\Longrightarrow $ (i) It is obvious (even if we have $m=1$).

(i) $\Longrightarrow $ (ii) Let us suppose that $m=\dim \mathcal{T}\geq 2$
and let us consider that $H$ is a Kronecker $h$-regular multi-time
Hamiltonian function, that is we have%
\begin{equation*}
{\frac{1}{2}}{\frac{\partial ^{2}H}{\partial p_{i}^{a}\partial p_{j}^{b}}}%
=h_{ab}(t^{c})g^{ij}(t^{c},x^{k},p_{k}^{c}).
\end{equation*}

($1^{\circ }$) Firstly, let us suppose that there exist two distinct indices 
$a$ and $b$, from the set $\{1,\ldots ,m\}$, such that $h_{ab}\neq 0$. Let $%
k $ ($c$, respectively) be an arbitrary element of the set $\{1,\ldots ,n\}$
($\{1,\ldots ,m\},$ respectively). Deriving the above relation, with respect
to the variable $p_{k}^{c}$, and using the Schwartz theorem, we obtain the
equalities%
\begin{equation*}
{\frac{\partial g^{ij}}{\partial p_{k}^{c}}}h_{ab}={\frac{\partial g^{jk}}{%
\partial p_{i}^{a}}}h_{bc}={\frac{\partial g^{ik}}{\partial p_{j}^{b}}}%
h_{ac},\quad \forall \;a,b,c\in \{1,\ldots ,m\},\quad \forall \;i,j,k\in
\{1,\ldots ,n\}.
\end{equation*}%
Contracting now with $h^{cd}$, we deduce that%
\begin{equation*}
{\frac{\partial g^{ij}}{\partial p_{k}^{c}}}h_{ab}h^{cd}=0,\quad \forall
\;d\in \{1,\ldots ,m\}.
\end{equation*}

In this context, the supposing $h_{ab}\neq 0$, together with the fact that
the metric $h$ is non-degenerate, imply that 
\begin{equation*}
{{\frac{\partial g_{ij}}{\partial p_{k}^{c}}}=0,}
\end{equation*}%
for any two arbitrary indices $k$ and $c$. Consequently, we have $%
g^{ij}=g^{ij}(t^{d},x^{r})$.

($2^{\circ }$) Let us suppose now that $h_{ab}=0,\;\forall \;a\neq b\in
\{1,\ldots ,m\}$. It follows that 
\begin{equation*}
h_{ab}=h_{a}(t)\delta _{b}^{a},\;\text{\ \ }\forall \;a,b\in \{1,\ldots ,m\},
\end{equation*}%
where $h_{a}(t)\neq 0,$\ $\forall \;a\in \{1,\ldots ,m\}.$ In these
conditions, the relations%
\begin{equation*}
\begin{array}{l}
{{\dfrac{\partial ^{2}L}{\partial p_{i}^{a}\partial p_{j}^{b}}}=0,\quad
\forall \;a\neq b\in \{1,\ldots ,m\},\quad \forall \;i,j\in \{1,\ldots ,n\}}%
,\medskip \\ 
{{\dfrac{1}{2h_{a}(t)}}{\dfrac{\partial ^{2}L}{\partial p_{i}^{a}\partial
p_{j}^{a}}}=g_{ij}(t^{c},x^{k},p_{k}^{c}),\quad \forall \;a\in \{1,\ldots
,m\},\quad \forall \;i,j\in \{1,\ldots ,n\},}%
\end{array}%
\end{equation*}%
are true. If we fix now an index $a$ in the set $\{1,\ldots ,m\}$, we deduce
from the first relations that the local functions ${{\dfrac{\partial L}{%
\partial p_{i}^{a}}}}$ depend only by the coordinates $%
(t^{c},x^{k},p_{k}^{a})$. Considering $b\neq a$ another index from the set $%
\{1,\ldots ,m\}$, the second relations imply%
\begin{equation*}
{{\dfrac{1}{2h_{a}(t)}}{\dfrac{\partial ^{2}L}{\partial p_{i}^{a}\partial
p_{j}^{a}}}={\dfrac{1}{2h_{b}(t)}}{\dfrac{\partial ^{2}L}{\partial
p_{i}^{b}\partial p_{j}^{b}}}=g_{ij}(t^{c},x^{k},p_{k}^{c}),\quad \forall
\;i,j\in \{1,\ldots ,n\}}.
\end{equation*}%
Because the first term of the above equality depends only by the coordinates 
$(t^{c},x^{k},p_{k}^{a})$, while the second term depends only by the
coordinates $(t^{c},x^{k},p_{k}^{b})$, and because we have $a\neq b$, we
conclude that $g^{ij}=g^{ij}(t^{d},x^{r})$.

Finally, the equalities%
\begin{equation*}
{\frac{1}{2}}{\frac{\partial ^{2}H}{\partial p_{i}^{a}\partial p_{j}^{b}}}%
=h_{ab}(t^{c})g^{ij}(t^{c},x^{k}),\quad \forall \;a,b\in \{1,\ldots
,m\},\quad \forall \;i,j\in \{1,\ldots ,n\},
\end{equation*}%
imply without difficulties that the multi-time Hamilton function $H$ is one
of non-autonomous electrodynamic kind (\ref{NEDTH}).
\end{proof}

\begin{corollary}
The \textit{fundamental vertical metrical d-tensor of a }Kronecker $h$%
-regular multi-time Hamiltonian function $H$ has the form%
\begin{equation}
G_{(a)(b)}^{(i)(j)}={\frac{1}{2}}{\frac{\partial ^{2}H}{\partial
p_{i}^{a}\partial p_{j}^{b}}}=\left\{ 
\begin{array}{ll}
h_{11}(t)g^{ij}(t,x^{k},p_{k}), & m=\dim \mathcal{T}=1\medskip \\ 
h_{ab}(t^{c})g^{ij}(t^{c},x^{k}), & m=\dim \mathcal{T}\geq 2.%
\end{array}%
\right.  \label{FVDT}
\end{equation}
\end{corollary}

\begin{remark}
i) It is obvious that the Theorem \ref{thchar} is an obstruction in the
development of a fertile geometrical theory of the multi-time Hamilton
spaces. This obstruction will be surpassed in other paper by the
introduction of the more general geometrical concept of \textbf{generalized
multi-time Hamilton space}. The generalized multi-time Hamilton geometry on
the dual 1-jet space $E^{\ast }$ will be constructed using only a Kronecker $%
h$-regular \textit{fundamental vertical metrical d-tensor (not necessarily
provided by a Hamiltonian function)}%
\begin{equation*}
{G_{(a)(b)}^{(i)(j)}=}h_{ab}(t^{c})g^{ij}(t^{c},x^{k},p_{k}^{c}),{\ }
\end{equation*}%
{together with an \textbf{a priori} given nonlinear connection }$N$ on $%
E^{\ast }${.}

ii) In the case $m=\dim \mathcal{T}\geq 2$, the Theorem \ref{thchar} obliges
us to continue our geometrical study of the multi-time Hamilton spaces
channeling our attention upon the \textbf{non-autonomous multi-time Hamilton
spaces of electrodynamics}.
\end{remark}

In the sequel, following the geometrical ideas of Miron from [9], we will
show that any Kronecker $h$-regular multi-time Hamiltonian function $H$
produces a natural nonlinear connection on the dual 1-jet bundle $E^{\ast }$%
, which depends only by $H$. In order to do that, let us take a Kronecker $h$%
-regular multi-time Hamiltonian function $H$, whose fundamental vertical
metrical d--tensor is given by (\ref{FVDT}). Also, let us consider the 
\textit{generalized spatial Christoffel symbols} of the d-tensor $g_{ij}$,
given by%
\begin{equation*}
\Gamma _{ij}^{k}=\frac{g^{kl}}{2}\left( \frac{\partial g_{li}}{\partial x^{j}%
}+\frac{\partial g_{lj}}{\partial x^{i}}-\frac{\partial g_{ij}}{\partial
x^{l}}\right) .
\end{equation*}

In this context, using preceding notations, we can give the following result:

\begin{theorem}
The pair of local functions $N=\left( \underset{1}{N}\text{{}}_{(i)b}^{(a)},%
\underset{2}{N}\text{{}}_{(i)j}^{(a)}\right) $ on $E^{\ast }$, where%
\begin{equation}
\begin{array}{l}
\underset{1}{N}\text{{}}_{(i)b}^{(a)}=\varkappa _{cb}^{a}p_{i}^{c},\medskip
\\ 
\underset{2}{N}\text{{}}_{(i)j}^{(a)}=\dfrac{h^{ab}}{4}\left[ \dfrac{%
\partial g_{ij}}{\partial x^{k}}\dfrac{\partial H}{\partial p_{k}^{b}}-%
\dfrac{\partial g_{ij}}{\partial p_{k}^{b}}\dfrac{\partial H}{\partial x^{k}}%
+g_{ik}\dfrac{\partial ^{2}H}{\partial x^{j}\partial p_{k}^{b}}+g_{jk}\dfrac{%
\partial ^{2}H}{\partial x^{i}\partial p_{k}^{b}}\right] ,%
\end{array}
\label{CNLC}
\end{equation}%
represents a nonlinear connection on $E^{\ast }$, which is called the 
\textbf{canonical nonlinear connection of the multi-time Hamilton space }$%
MH_{m}^{n}=(E^{\ast },H).$
\end{theorem}

\begin{proof}
Taking into account the classical transformation rules of the Christoffel
symbols $\varkappa _{bc}^{a}$ of the temporal semi-Riemannian metric $%
h_{ab}, $ by direct local computations, we deduce that the temporal
components $\underset{1}{N}${}$_{(i)b}^{(a)}$ from (\ref{CNLC}) verify the
first transformation rules from (\ref{schco}) (please see also the Corollary %
\ref{NLC0}).

In the particular case when $m=\dim \mathcal{T}=1,$ the spatial components%
\begin{equation*}
\underset{2}{N}\text{{}}_{(i)j}^{(1)}=\dfrac{h^{11}}{4}\left[ \dfrac{%
\partial g_{ij}}{\partial x^{k}}\dfrac{\partial H}{\partial p_{k}}-\dfrac{%
\partial g_{ij}}{\partial p_{k}}\dfrac{\partial H}{\partial x^{k}}+g_{ik}%
\dfrac{\partial ^{2}H}{\partial x^{j}\partial p_{k}}+g_{jk}\dfrac{\partial
^{2}H}{\partial x^{i}\partial p_{k}}\right]
\end{equation*}%
become (except the multiplication factor $h^{11}$) exactly the canonical
nonlinear connection from the classical Hamilton geometry (please see [9] or
[11, pp. 127]).

For $m=\dim \mathcal{T}\geq 2,$ the Theorem \ref{thchar} (more exactly, the
formula (\ref{NEDTH})) leads us to the following expression for the spatial
components $\underset{2}{N}${}$_{(i)j}^{(a)}$ from (\ref{CNLC}):%
\begin{equation}
\underset{2}{N}\text{{}}_{(i)j}^{(a)}=-\Gamma
_{ij}^{k}p_{k}^{a}+T_{(i)j}^{(a)},  \label{ExprCNLC}
\end{equation}%
where%
\begin{equation*}
T_{(i)j}^{(a)}=\dfrac{h^{ab}}{4}\left[ \dfrac{\partial g_{ij}}{\partial x^{k}%
}U_{(b)}^{(k)}+g_{ik}\dfrac{\partial U_{(b)}^{(k)}}{\partial x^{j}}+g_{jk}%
\dfrac{\partial U_{(b)}^{(k)}}{\partial x^{i}}\right] .
\end{equation*}

Because $T_{(i)j}^{(a)}$ is a d-tensor on $E^{\ast }$ (we prove this by
local computations, studying the transformation laws of $T_{(i)j}^{(a)}$),
it immediately follows that the spatial components $\underset{2}{N}${}$%
_{(i)j}^{(a)}$ given by (\ref{ExprCNLC}) transform as in the second laws of (%
\ref{schco}).
\end{proof}

Finally, using the expression (\ref{ExprCNLC}), by computations, we find

\begin{corollary}
For $m=\dim \mathcal{T}\geq 2,$ the canonical nonlinear connection $N$ of a
multi-time Hamilton space $MH_{m}^{n}=(E^{\ast },H)$ (given by (\ref{NEDTH}%
)) has the components%
\begin{equation*}
\begin{array}{l}
\underset{1}{N}\text{{}}_{(i)b}^{(a)}=\varkappa _{cb}^{a}p_{i}^{c},\medskip
\\ 
\underset{2}{N}\text{{}}_{(i)j}^{(a)}=-\Gamma _{ij}^{k}p_{k}^{a}+\dfrac{%
h^{ab}}{4}\left( U_{ib\bullet j}+U_{jb\mathbf{\bullet }i}\right) ,%
\end{array}%
\end{equation*}%
where $U_{ib}=g_{ik}U_{(b)}^{(k)}$ and 
\begin{equation*}
U_{kb\mathbf{\bullet }r}=\dfrac{\partial U_{kb}}{\partial x^{r}}%
-U_{sb}\Gamma _{kr}^{s}.
\end{equation*}
\end{corollary}

\textbf{Author's address: }Gheorghe ATANASIU, Str. Gh. Baiulescu, Nr. 11, Bra%
\c{s}ov, BV 500107, Romania.

\textbf{E-mail}: gh\_atanasiu@yahoo.com

\textbf{Place of work: }University "Transilvania" of Bra\c{s}ov, Faculty of
Mathematics and Informatics.

\bigskip

\textbf{Author's address: }Mircea NEAGU, Str. L\u{a}m\^{a}i\c{t}ei, Nr. 66,
Bl. 93, Sc. G, Ap. 10, Bra\c{s}ov, BV 500371, Romania.

\textbf{E-mail}: mirceaneagu73@yahoo.com

\textbf{Place of work: }University "Transilvania" of Bra\c{s}ov, Faculty of
Mathematics and Informatics.

\end{document}